\documentclass[10pt,leqno]{amsart}
\usepackage{graphicx}
\usepackage{indentfirst,csquotes}

\usepackage{lipsum} 
\usepackage{amssymb,amsthm,amsmath}
\usepackage{xcolor,paralist,hyperref,titlesec,fancyhdr,etoolbox}
\newtheorem{theorem}{Theorem}[]

\newtheorem{remark}[theorem]{Remark}

\titleformat{\section}[display]{\normalfont\huge\bfseries\centering}{\centering\chaptertitlename\thechapter}{10pt}{\Large}
\titlespacing*{\section}{0pt}{0ex}{0ex}

\hypersetup{ colorlinks=true, linkcolor=black, filecolor=black, urlcolor=black }

\usepackage{lipsum}

\begin{document}

\title[LIOUVILLE TYPE THEOREMS FOR ISENTROPIC NAVIER-STOKES SYSTEM]{ON LIOUVILLE TYPE THEOREMS FOR THE 3D STEADY ISENTROPIC NAVIER-STOKES SYSTEM WITHOUT D-CONDITION} 

\author[Jiu]{Quansen Jiu}
\date{\today}
\address{Beijing, People's Republic of China}
\email{jiuqs@cnu.edu.cn}

\author[Tan]{Jie Tan}
\date{\today}
\address{Beijing, People's Republic of China}
\email{2250501015@cnu.edu.cn}

\author[Yan]{Zhihong Yan}
\date{\today}
\address{Beijing, People's Republic of China}
\email{yzh00620math@163.com}

\let\thefootnote\relax

\begin{abstract}
In this paper, we establish Liouville-type theorems for the steady compressible Navier-Stokes system. Assuming a smooth solution \(u \in L^p(\mathbb{R}^3)\), \(3 \le p \le \frac{9}{2}\), with bounded density, one obtains \(u \equiv0\). This generalizes the result of Li-Yu \cite{Li-Yu} by removing the Dirichlet condition \(\int_{\mathbb{R}^3} |\nabla u|^2 \, dx < \infty\). If \(\frac{9}{2} < p < 6\),  Liouville-type theorem holds under the additional oscillation condition for momentum \(\rho u \in \dot{B}^{\frac{3}{p} - \frac{3}{2}}_{\infty,\infty}(\mathbb{R}^3)\). For the marginal case \(u \in L^6(\mathbb{R}^3)\), the oscillation condition can be replaced by \(\rho u \in BMO^{-1}(\mathbb{R}^3)\). We also present results in Morrey-type spaces: \(u \in \dot{M}^{s,6}(\mathbb{R}^3)\) and \(\rho u \in \dot{M}_w^{q,3}(\mathbb{R}^3)\) for \(2 \le s \le 6\) and \(\frac{3}{2} < q \le 3\). 
\end{abstract} 

\bigskip
\maketitle
\section*{Introduction}

The three dimensional steady isentropic compressible Navier-Stokes equations are
 \begin{equation}\label{eqNS}
 	\begin{cases}
\mathrm{div}(\rho u \otimes u) + \nabla {P}(\rho ) = \mu\Delta  u+(\mu+\lambda)\nabla \mathrm{div}u,\\[1ex]
\mathrm{div}(\rho u) = 0,
 	\end{cases} 		
 \end{equation}
where $\rho  \ge 0$ is density, $u=(u_1,u_2,u_3)$ stands for velocity, constants $\mu> 0$
and $\mu+\lambda> 0$ are viscosity coefficients, and the pressure-density function takes form
\begin{align}\label{dfqt}
	{p}(\rho ) = {A}{\rho ^\gamma }
\end{align}
with ${A} > 0$, $\gamma  > 1$.

Significant progress has recently been made regarding Liouville-type theorems for the steady incompressible Navier--Stokes system.
A classical result in this direction, presented in the monograph by Galdi~\cite{Galdi2011}, establishes that any smooth solution
must satisfy $u \equiv 0$ under the decay condition $u \in L^{\frac{9}{2}}(\mathbb{R}^3)$.
The same conclusion remains valid under the weaker integrability condition $u \in L^p(\mathbb{R}^3)$ for $3 \le p < \frac{9}{2}$,
as shown by Chamorro et al.~\cite{D.Chamorro.2021}.
An important refinement was introduced by Seregin~\cite{seregin L^6 BMO^-1}, where the decay condition is replaced by a
decay--oscillation condition: $u \in L^6(\mathbb{R}^3) \cap BMO^{-1}(\mathbb{R}^3)$. This result has subsequently been extended:
Chamorro et al.~\cite{D.Chamorro.2021} introduced a new decay--oscillation criterion of the form
\[
u \in L^p(\mathbb{R}^3) \cap \dot{B}_{\infty,\infty}^{\frac{3}{p}-\frac{3}{2}}(\mathbb{R}^3), \quad \frac{9}{2} < p < 6.
\]
More recently, Seregin~\cite{Seregin-16} obtained a Liouville-type theorem in Morrey-type spaces, specifically under the conditions
\[
u \in \dot{M}^{s,6}(\mathbb{R}^3) \cap \dot{M}_w^{q,3}(\mathbb{R}^3).
\]


For the compressible Navier--Stokes system, the literature on Liouville-type theorems remains relatively sparse---see, e.g., \cite{Chae, Li-Niu, Li-Yu}.
In the setting of smooth $D$-solutions (i.e., solutions satisfying $u \in \dot{H}^1(\mathbb{R}^3)$), D. Chae \cite{Chae} established a Liouville-type theorem under the integrability condition $u \in L^{\frac{3}{2}}(\mathbb{R}^3)$. Subsequently, D. Li et al. \cite{Li-Yu} proved that any smooth $D$-solution must vanish identically provided $u \in L^{\frac{9}{2}}(\mathbb{R}^3)$. This $L^{\frac{9}{2}}$ condition was later relaxed by Z.-Y. Li et al. \cite{Li-Niu}, who replaced it with a suitable Lorentz-space assumption.
To our knowledge, however, no previous results have removed the Dirichlet-type regularity condition $u \in \dot{H}^1(\mathbb{R}^3)$ entirely. The main objective of the present paper is to investigate Liouville-type theorems for the steady compressible Navier--Stokes system \emph{without} imposing such a $D$-condition.

In our first theorem, we will prove Liouville type properties for the system \eqref{eqNS}-\eqref{dfqt} under decay condition.
\begin{theorem}\label{thNS1}
	Let $ (u,\rho) $ is a smooth solution of \eqref{eqNS}-\eqref{dfqt} which satisfies $ u\in L^p(\mathbb{R}^3)$ and $\rho\in L^{\infty}(\mathbb{R}^3)$ with $3\leq p\leq\frac92$, then $u=0$ and $\rho=$constant.
\end{theorem}

\begin{remark}
       For the incompressible case, Galdi's result~\cite{Galdi2011} established the Liouville theorem under the condition $u \in L^{\frac{9}{2}}(\mathbb{R}^3)$, which can be viewed as a limiting function space for this problem. This result was later extended by Chamorro et al., who proved the theorem for the larger range $u \in L^p(\mathbb{R}^3)$ with $3 \le p \le \frac{9}{2}$.
In the compressible setting, Chae~\cite{Chae} obtained a Liouville theorem for $D$-solutions satisfying $u \in L^{\frac{3}{2}}(\mathbb{R}^3)$. Subsequently, Li et al.~\cite{Li-Yu} relaxed this condition to $u \in L^{\frac{9}{2}}(\mathbb{R}^3)$. The present work improves on the result of Li et al.~\cite{Li-Yu} by completely removing the $D$-condition.
\end{remark}

A crucial observation is that for a vector field \( b \) satisfying \( \operatorname{div} b = 0 \), one may represent it in the form \( b = \operatorname{div} d \) for some antisymmetric matrix \( d = (d_{ij}) \). More details on this representation can be found in the work of Seregin et al.~\cite{Seregin-Silvestre-Sverak-Zlatos}.
Based on this fact, we say that a divergence-free vector field \( b \) belongs to \( BMO^{-1}(\mathbb{R}^3) \) if there exists an antisymmetric matrix \( d \in BMO(\mathbb{R}^3) \) such that \( b = \operatorname{div} d \).
Additionally, the homogeneous Besov space of the negative index can be characterized as in Bahouri et al.~\cite{Bahouri}:
\[
\dot{B}_{\infty,\infty}^q(\mathbb{R}^3)=\Bigl\{ f \in \mathcal{S}' \;:\; \| f \|_{\dot{B}_{\infty,\infty}^q}
= \sup_{t>0} \, t^{-q/2} \| h_t * f \|_{L^\infty} < \infty \Bigr\}, \quad q < 0,
\]
where \( h_t \) denotes the heat kernel at time \( t \).

\begin{theorem}\label{thNS2}
	Assume that $ (u,\rho) $ is a smooth solution to the system \eqref{eqNS}-\eqref{dfqt} in $\mathbb{R}^3$.
	\begin{enumerate}
		\item[(1)] If $ u\in L^6(\mathbb{R}^3)$ and $\rho u\in BMO^{-1}(\mathbb{R}^3) $, then $u=0$ and $\rho=$constant.
		\item[(2)] If $ u\in L^p(\mathbb{R}^3)$ and $\rho u\in \dot{B}_{\infty,\infty}^{\frac3p-\frac32}(\mathbb{R}^3)$ with $\frac92<p<6$, then $u=0$ and $\rho=$constant.
	\end{enumerate}
\end{theorem}

\begin{remark}
       For the incompressible case, Seregin~\cite{seregin L^6 BMO^-1} first established a Liouville theorem under the endpoint condition
\( u \in L^6(\mathbb{R}^3) \cap BMO^{-1}(\mathbb{R}^3) \). This result was later extended to the larger class
\[
u \in L^p(\mathbb{R}^3) \cap \dot{B}_{\infty,\infty}^{\frac{3}{p} - \frac{3}{2}}(\mathbb{R}^3), \quad 3 \le p \le \frac{9}{2}.
\]
In the compressible setting, the vector field \( \rho u \) is divergence-free. Consequently, it is natural to impose the oscillation conditions directly on the momentum \( \rho u \) rather than on the velocity field alone.

\end{remark}

For \(1 < p \le q < \infty\), define the homogeneous Morrey space that
\[
\dot{M}^{p,q}(\mathbb{R}^3):=
\Bigl\{ f \in L^p_{\mathrm{loc}}(\mathbb{R}^3):\;
\|f\|_{\dot{M}^{p,q}(\mathbb{R}^3)} =
\sup_{x \in \mathbb{R}^3,\;R>0} R^{3\left(\frac{1}{q}-\frac{1}{p}\right)}
\|f\|_{L^{p}(B_R(x))} < \infty \Bigr\},
\]
where \(B_R(x)\) denotes the ball of radius \(R\) centered at \(x\).
If the strong \(L^p\)-norm in the definition is replaced by the weak \(L^p\)-norm
\[
\|f\|_{L_w^{p}(B_R(x))} :=
\sup_{t>0} t \, \bigl|\{ y \in B_R(x): |f(y)| > t \}\bigr|^{1/p},
\]
one obtains the weak Morrey space
\[
\dot{M}_w^{p,q}(\mathbb{R}^3) :=
\Bigl\{ f \in L^{p}_{w,\mathrm{loc}}(\mathbb{R}^3):\;
\|f\|_{\dot{M}_w^{p,q}(\mathbb{R}^3)} =
\sup_{x \in \mathbb{R}^3,\;R>0} R^{3\left(\frac{1}{q}-\frac{1}{p}\right)}
\|f\|_{L_w^{p}(B_R(x))} < \infty \Bigr\}.
\]
Recall that the weak Lebesgue space \(L^{p,\infty}(\Omega)\) is a special case of the Lorentz space \(L^{p,m}(\Omega)\) (with \(m=\infty\)), whose norm is given by
\[
\|f\|_{L^{p,m}(\Omega)} =
\begin{cases}
\displaystyle
\Bigl(\int_0^{+\infty} t^{m-1} \bigl|\{y\in\Omega:|f(y)|>t\}\bigr|^{\frac{m}{p}} \, dt \Bigr)^{\frac{1}{m}}, & \text{if } 1 \le m < +\infty, \\[10pt]
\displaystyle
\sup_{t>0} \, t \, \bigl|\{y\in\Omega:|f(y)|>t\}\bigr|^{\frac{1}{p}}, & \text{if } m = +\infty.
\end{cases}
\]
We now state our main result.

\begin{theorem}\label{thNS3}
Let \((u,\rho)\) be a smooth solution of the steady compressible Navier--Stokes system in \(\mathbb{R}^3\). Assume that
$
u \in \dot{M}^{s,6}(\mathbb{R}^3)\; \text{and}\;
\rho u \in \dot{M}_w^{q,3}(\mathbb{R}^3)$,
with parameters that satisfy
$2 \le s \le 6, \; \frac{3}{2} < q \le 3.$
Then \(u \equiv 0\) and \(\rho\) is constant throughout \(\mathbb{R}^3\).
\end{theorem}


\begin{remark}
  Due to properties of Morrey-type space, we see $L^6(\mathbb{R}^3)$ is substituted by the Morrey space $\dot{M}^{s,6}(\mathbb{R}^3)$ and the space $\dot M_w^{q,3} ( {\mathbb{R}^3} )$
  has the same degree with $BMO^{-1}(\mathbb{R}^3)$.
\end{remark}

\section*{Proofs of Main Theorems }

\begin{proof}[{\rm \textbf{Proof of Theorem \ref{thNS1}}}]
	
Without loss of generality we assume that $ A=1. $ Fix $ R>0 $ and choose $ \frac R2\leq\rho<r\leq R $. Let a radial cut-off function $ \phi\in C_0^{\infty}(\mathbb{R}^3) $ possess the properties as follows:
\begin{equation*}
  0\leq\phi\leq1,\quad\phi=1\ \text{in}\ B_{\rho},\quad\phi=1\ \text{in}\ B_r^c,\quad |\nabla\phi|\leq\frac{c}{r-\rho},\quad|\nabla^2\phi|\leq\frac{c}{(r-\rho)^2},
\end{equation*}
where $ B_r$ denote the ball $ B(0,r)\subset\mathbb{R}^3 $.

For any $ s>1 $, there exists a constant $ c_0=c_0(s,\gamma)>0 $ and a function $w\in W_0^{1,s}(B_r)$,  such that $\mathrm{div}w=\frac{\gamma }{{\gamma  - 1}}\nabla\phi\cdot u$ and
\begin{equation}\label{eq2.th1.1}
	\int_{B_r}|\nabla w|^s\mathrm{d}x\leq c_0\int_{B_r}|\nabla\phi\cdot u|^s\mathrm{d}x\leq\frac{c_0}{(r-\rho)^s}\int_{B_{R}\setminus B_{R/2}}| u|^s\mathrm{d}x.
\end{equation}
Multiplying $\eqref{eqNS}_1$ by $\phi u-w$ and integrating the resultant over $B_r$ yields
\begin{equation}\label{eq2.th1.2}
	\begin{aligned}
		\mu\int_{B_r}\phi|\nabla u|^2&\mathrm{d}x+(\mu+\lambda)\int_{B_r}\phi|\mathrm{div}u|^2\mathrm{d}x\\
		&=\frac12\int_{B_r}\rho|u|^2 u\cdot\nabla\phi\mathrm{d}x+\int_{B_r}\rho (u\cdot\nabla u)\cdot w\mathrm{d}x\\
        &\quad+\int_{B_r}(\nabla\rho^{\gamma}\cdot w-\gamma\phi\rho^{\gamma-1}u\cdot\nabla\rho)\mathrm{d}x-\mu\int_{B_r}\nabla\phi\cdot(\nabla u\cdot u)\mathrm{d}x\\
        &\quad-\mu\int_{B_r}\nabla u:\nabla w\mathrm{d}x+(\mu+\lambda)\int_{B_r}\mathrm{div}u(\mathrm{div}w-u\cdot\nabla\phi)\mathrm{d}x.\\
	\end{aligned}
\end{equation}
Integrating by parts, it holds that
\begin{equation}\label{ttwo}
  \int_{B_r}\rho (u\cdot\nabla u)\cdot w\mathrm{d}x=-\int_{B_r} \rho u\cdot(\nabla w\cdot u)\mathrm{d}x.
\end{equation}
Thanks to the identity $\mathrm{div}w=\frac{\gamma }{{\gamma  - 1}}\nabla\phi\cdot u$, one obtains
\begin{equation}\label{tthree}
\begin{aligned}
	\int_{B_r}\nabla\rho^{\gamma}\cdot w\mathrm{d}x-\gamma\int_{B_r}&\phi\rho^{\gamma-1}u\cdot\nabla\rho\mathrm{d}x\\
&=-\frac{\gamma}{\gamma-1}\int_{B_r}\phi\rho u\cdot\nabla(\rho^{\gamma-1})\mathrm{d}x-\int_{B_r}\rho^{\gamma}\mathrm{div}w\mathrm{d}x\\
	&=\frac{\gamma}{\gamma-1}\int_{B_r}\rho^{\gamma} u\cdot\nabla\phi\mathrm{d}x-\int_{B_r}\rho^{\gamma}\mathrm{div}w\mathrm{d}x=0.
\end{aligned}
\end{equation}
By straightforward calculations, it follows from property of $w$ that
\begin{align}\label{tfour}
  (\mu+\lambda)\int_{B_r}\mathrm{div}u(\mathrm{div}w-u\cdot\nabla\phi)\mathrm{d}x=-\frac{(\mu+\lambda)}{\gamma-1}\int_{B_r}\mathrm{div}u\cdot  (u\cdot\nabla\phi)\mathrm{d}x.
\end{align}
Plugging \eqref{ttwo}-\eqref{tfour} into \eqref{eq2.th1.2}, we conclude
\begin{equation}\label{eq2.th1.3}
	\begin{aligned}
		\mu\int_{B_r}\phi|\nabla u|^2&\mathrm{d}x+(\mu+\lambda)\int_{B_r}\phi|\mathrm{div}u|^2\mathrm{d}x\\
		&=\frac12\int_{B_r}\rho|u|^2 u\cdot\nabla\phi\mathrm{d}x-\int_{B_r} \rho u\cdot(\nabla w\cdot u)\mathrm{d}x\\
        &\quad-\mu\int_{B_r}\nabla\phi\cdot(\nabla u\cdot u)\mathrm{d}x-\mu\int_{B_r}\nabla u:\nabla w\mathrm{d}x\\
        &\quad-\frac{(\mu+\lambda)}{\gamma-1}\int_{B_r}\mathrm{div}u\cdot u\cdot\nabla\phi\mathrm{d}x = I_1+I_2+I_3+I_4+I_5.
	\end{aligned}
\end{equation}
It is clear that
\begin{equation}\label{I1}
  |I_1|\leq \frac{c\Vert\rho\Vert_{L^{\infty}}}{r-\rho}\Vert u\Vert_{L^3(B_{R}\setminus B_{R/2})}^3.
\end{equation}
Applying the inequality \eqref{eq2.th1.1} for the case $s=3$ yields
\begin{equation}\label{I2}
  |I_2|\leq c\Vert\rho\Vert_{L^{\infty}}\int_{B_r}|u|^2|\nabla w|\mathrm{d}x\leq \frac{c\Vert\rho\Vert_{L^{\infty}}}{r-\rho}\Vert u\Vert_{L^3(B_R)}^2\Vert u\Vert_{L^3(B_{R}\setminus B_{R/2})},
\end{equation}
and
\begin{equation}\label{I3}
  |I_3|+|I_4|\leq\frac{\mu}{4}\int_{B_r}|\nabla u|^2\mathrm{d}x+\frac{c}{(r-\rho)^2}\Vert u\Vert_{L^2(B_{r}\setminus B_{\rho})}^2.
\end{equation}
It follows from the Young's inequality that
\begin{equation}\label{I5}
  |I_5|\leq\frac{\mu+\lambda}{4}\int_{B_r} |\mathrm{div}u|^2\mathrm{d}x+\frac{c}{(r-\rho)^2}\Vert u\Vert_{L^2(B_{r}\setminus B_{\rho})}^2.
\end{equation}
Combining \eqref{eq2.th1.3}-\eqref{I5} and using H\"older inequality yield
\begin{equation}\label{eq2.th1.4}
	\begin{aligned}
		\int_{B_{\rho}}(\mu|\nabla u|^2+(\mu+\lambda)|\mathrm{div}u|^2)\mathrm{d}x
		\leq&\frac14\int_{B_r}(\mu|\nabla u|^2+(\mu+\lambda)|\mathrm{div}u|^2)\mathrm{d}x\\
		&+\frac{c}{r-\rho}\left(\Vert u\Vert_{L^3(B_{R}\setminus B_{R/2})}^3+\Vert u\Vert_{L^3(B_R)}^2\Vert u\Vert_{L^3(B_{R}\setminus B_{R/2})}\right)\\
        &+\frac{cR}{(r-\rho)^2}\Vert u\Vert_{L^3(B_{R}\setminus B_{R/2})}^2,
	\end{aligned}
\end{equation}
where $c=c(c_0, \mu, \lambda, \Vert\rho\Vert_{L^{\infty}})$.
For convenience, we denote
\begin{gather*}
	\Phi_u(r)=\int_{B_r}(\mu|\nabla u|^2+(\mu+\lambda)|\mathrm{div}u|^2)\mathrm{d}x,
\end{gather*}
and
\begin{equation*}
  \kappa_R^1= c\left(\Vert u\Vert_{L^3(B_{R}\setminus B_{R/2})}^3+\Vert u\Vert_{L^3(B_R)}^2\Vert u\Vert_{L^3(B_{R}\setminus B_{R/2})}\right), \,\,\,
    \kappa_R^2= cR\Vert u\Vert_{L^3(B_{R}\setminus B_{R/2})}^2.
\end{equation*}
The energy inequality \eqref{eq2.th1.4} arrives at
\begin{equation}\label{eq2.th1.5}
	\Phi_u(\rho)\leq\frac14\Phi_u(r)+\frac{\kappa_R^1}{r-\rho}+\frac{\kappa_R^2}{(r-\rho)^2},\quad\text{for any } \frac{R}{2}\leq\rho<r\leq R.
\end{equation}
For $ k\in\mathbb{N}, $ letting $ r_k= R2^{-1/k} $ and substituting it into inequality \eqref{eq2.th1.5}, one has
\begin{align*}
	\Phi_u(r_k)&\leq\frac14\Phi_u(r_{k+1})+\kappa_R^1R^{-1}(2^{-1/(k+1)}-2^{-1/k})^{-1}+\kappa_R^2R^{-2}(2^{-1/(k+1)}-2^{-1/k})^{-2}\\
	&\leq\frac14\Phi_u(r_{k+1})+c\kappa_R^1R^{-1}k^2+c\kappa_R^2R^{-2}k^{4},
\end{align*}
where we used the fact that
\begin{equation*}
  2^{-1/(k+1)}-2^{-1/k}=2^{-1/k}(2^{1/k(k+1)}-1)\geq\frac{c\log2}{2k^2}.
\end{equation*}
Thanks to iteration argument, it holds that
\begin{align}\label{iain}
	\Phi_u(r_1)&\leq\frac1{4^k}\Phi_u(r_{k+1})+c\kappa_R^1R^{-1}\sum_{j=1}^{k}\frac{j^2}{4^{j-1}}+c\kappa_R^2R^{-2}\sum_{j=1}^{k}\frac{j^{4}}{4^{j-1}}.
\end{align}
Passing \eqref{iain} as $k \to \infty$, since $\sum_{j=1}^{\infty}\frac{j^2}{4^{j-1}}<+\infty$ and $\sum_{j=1}^{\infty}\frac{j^{4}}{4^{j-1}}<\infty$, Caccioppoli type inequality can be obtained
\begin{equation}\label{eq2.th1.6}
	\begin{aligned}
		\int_{B_{R/2}}(\mu|\nabla u|^2+&(\mu+\lambda)|\mathrm{div}u|^2)\mathrm{d}x\\
		&\leq \frac cR\left(\Vert u\Vert_{L^3(B_{R}\setminus B_{R/2})}^3+\Vert u\Vert_{L^3(B_R)}^2\Vert u\Vert_{L^3(B_{R}\setminus B_{R/2})}+\Vert u\Vert_{L^3(B_{R}\setminus B_{R/2})}^2\right).
	\end{aligned}
\end{equation}
For $p=3$, the right hand side of inequality \eqref{eq2.th1.6} tends to zero as $R\to\infty$, then $u=0$.
If $3< p\leq\frac92$, one deduces
\begin{equation}\label{eq2.th1.7}
  \begin{aligned}
	\int_{B_{R/2}}(\mu|&\nabla u|^2+(\mu+\lambda)|\mathrm{div}u|^2)\mathrm{d}x\\
	&\leq  cR^{3\frac{p-3}{p}-1}\left(\Vert u\Vert_{L^p(B_{R}\setminus B_{R/2})}^3+\Vert u\Vert_{L^p(B_R)}^2\Vert u\Vert_{L^p(B_{R}\setminus B_{R/2})}\right)\\
    &\quad+cR^{2\frac{p-3}{p}-1}\Vert u\Vert_{L^p(B_{R}\setminus B_{R/2})}^2.
\end{aligned}
\end{equation}
Note that $ 3\frac{p-3}{p}-1\leq0$ and $2\frac{p-3}{p}-1<0 $. Passing limit to \eqref{eq2.th1.7} as $R\to\infty$ yields $u=0,\rho=\text{constant},$ so the conclusion follows.

\end{proof}

\begin{proof}[{\rm \textbf{Proof of Theorem \ref{thNS2}}}]
	\textbf{The proof of the first part.} Similar to the identity \eqref{eq2.th1.2}, one has
	\begin{equation}\label{eq2.th2.1}
		\begin{aligned}
			\mu\int_{B_r}\phi|\nabla u|^2&\mathrm{d}x+(\mu+\lambda)\int_{B_r}\phi|\mathrm{div}u|^2\mathrm{d}x\\
			&=-\mu\int_{B_r}\nabla u:(u\otimes\nabla \phi)\mathrm{d}x-\mu\int_{B_r}\nabla u:\nabla w\mathrm{d}x\\
			&\qquad-\frac{(\mu+\lambda)}{\gamma-1}\int_{B_r}\mathrm{div}u\cdot(u\cdot\nabla \phi)\mathrm{d}x-\int_{B_r}(\rho u\cdot\nabla u)\cdot(u\phi-w)\mathrm{d}x\\
			& = I_3+I_4+I_5+I_6.
		\end{aligned}
	\end{equation}
Due to $ \rho u\in BMO^{-1} $, there exists an anti-symmetric matrix $ d=(d_{ij})\in BMO $ such that
\begin{equation*}
  \rho u=\mathrm{div}d=(d_{ij,j}).
\end{equation*}

It is known in Stein \cite{Stein E M.1993} that if $d \in BMO$, then
\begin{equation*}
  \Gamma \left( s \right) =\sup_{x\in\mathbb{R}^3,r>0}\left(\frac1{|B(x,r)|}\int_{B(x,r)}|d-[d]_{B(x,r)}|^s\right)^{\frac1s}<\infty,
\end{equation*}
for each $1 \le s< \infty $, where $ [d]_{B(x,r)} $ is the mean value of $ d $ over $B(x,r)$.

For $s>2$, it holds that
\begin{equation}\label{eq2.th2.2}
	\begin{aligned}
		|I_6|&=\left|\int_{B_r}\bar{d}_{ij,j}u_{m,i}(u_m\phi-w_m)\mathrm{d}x\right|\\
		&\leq\int_{B_r}|\bar{d}_{ij}u_{m,i}u_{m}\phi_{,j}|\mathrm{d}x+\int_{B_r}|\bar{d}_{ij}u_{m,i}w_{m,j}|\mathrm{d}x\\
		&\leq cR^{\frac{3s-6}{2s}}\Vert \nabla u\Vert_{L^2({B_r})}(\Vert  \nabla w\Vert_{L^s({B_r})}+\frac{1}{r-\rho}\Vert  u\Vert_{L^s({B_R\setminus B_{R/2}})})\Gamma \left({\frac{2s}{s-2}} \right)\\
		&\leq\frac{cR^{\frac{3s-6}{2s}}}{r-\rho}\Vert \nabla u\Vert_{L^2({B_r})}\Vert  u\Vert_{L^s({B_R\setminus B_{R/2}})}\Gamma \left({\frac{2s}{s-2}} \right),
	\end{aligned}
\end{equation}
where the indices after comma mean derivatives and we denote $\bar{d}= d-[d]_{B_r}$.
For another terms, one obtains
\begin{gather}\label{eq2.th2.3}
|I_3+I_4+I_5|\leq \frac{cR^{\frac{3s-6}{2s}}}{r-\rho}\Vert  u\Vert_{L^s({B_R\setminus B_{R/2}})}\left(\Vert \nabla u\Vert_{L^2({B_r})}+\Vert \mathrm{div} u\Vert_{L^2({B_r})}\right).
\end{gather}
Combining (\ref{eq2.th2.1})-(\ref{eq2.th2.3}) yields
\begin{align*}
	\mu\int_{B_r}&\phi|\nabla u|^2\mathrm{d}x+(\mu+\lambda)\int_{B_r}\phi|\mathrm{div}u|^2\mathrm{d}x\\
 & \leq\frac{cR^{\frac{3s-6}{2s}}}{r-\rho}\Vert \mathrm{div} u\Vert_{L^2({B_r})}\Vert u\Vert_{L^s({B_R\setminus B_{R/2}})}
+\frac{c(s)R^{\frac{3s-6}{2s}}}{r-\rho}\Vert\nabla u\Vert_{L^2({B_r})}\Vert  u\Vert_{L^s({B_R\setminus B_{R/2}})}.
\end{align*}
Thanks to Young's inequality, it holds that
\begin{gather*}
	\Phi_u(\rho)\leq\frac14\Phi_u(r)+c(s)\frac{R^{\frac{3s-6}{s}}}{(r-\rho)^2}\Vert  u\Vert_{L^s({B_R\setminus B_{R/2}})}^2,
\end{gather*}
which implies the Caccioppoli type inequality by suitable iteration
\begin{equation*}
  \mu\int_{B_{R/2}}|\nabla u|^2\mathrm{d}x+(\mu+\lambda)\int_{B_{R/2}}|\mathrm{div}u|^2\mathrm{d}x\leq c(s)R^{\frac{s-6}{s}}\Vert  u\Vert_{L^s({B_R\setminus B_{R/2}})}^2.
\end{equation*}
Choosing $s=6$ and passing the limit as $R \to \infty$, we complete the proof of the first part.

\textbf{The proof of the second part.} The estimate of $|I_3+I_4+I_5|$ is as same as the proof of the first part. But the control of term $I_6$ is different.

Due to vector $\rho u$ is divergence-free, it holds that
\begin{equation*}
  \Delta(\rho u)=-\mathrm{curl}\,\mathrm{curl}(\rho u)+\nabla\mathrm{div}(\rho u)=-\mathrm{curl}\,\mathrm{curl}(\rho u).
\end{equation*}
Denoting $ v= (-\Delta)^{-1}\mathrm{curl}(\rho u)$, it is easy to check
\begin{equation*}
  \mathrm{curl}v=(-\Delta)^{-1}(\mathrm{curl}\,\mathrm{curl}(\rho u))=(-\Delta)^{-1}(-\Delta(\rho u))=\rho u.
\end{equation*}
Let $ v^*= v-v(0). $ Then $ \rho u_i=v^*_{j,k}-v^*_{k,j}(i=1,2,3) $ with $ (ijk)=(123)\in S_3 $, the symmetric group on 3 letters. Substituting it into $ I_6 $ and integrating by parts, by the appendix of Chamorro et al. \cite{D.Chamorro.2021}, we deduce
\begin{align*}
	I_6&=-\sum_{i,m=1}^{3}\sum_{(ijk)=(123)}\int_{B_r}(\rho u_i)u_{m,i}(\phi u_m-w_m)\mathrm{d}x\\
      &=-\sum_{i,m=1}^{3}\sum_{(ijk)=(123)}\int_{B_r}(v^*_{j,k}-v^*_{k,j})u_{m,i}(\phi u_m-w_m)\mathrm{d}x\\
	  &=\sum_{i,m=1}^{3}\sum_{(ijk)=(123)}\int_{B_r}\left(u_mu_{m,i}(v^*_{j}\phi _{,k}-v^*_{k}\phi _{,j})+u_{m,i}(v^*_{k}w_{m,j}-v^*_{j}w_{m,k})\right)\mathrm{d}x,
\end{align*}
which implies
\begin{equation}\label{eq2.th2.4}
	\begin{aligned}
		|I_6|&\leq c\left(\int_{B_r}|u||\nabla u||v^*||\nabla\phi|\mathrm{d}x+\int_{B_r}|\nabla u||v^*||\nabla w|\mathrm{d}x\right)\\
		&\leq \frac{c}{r-\rho}\left(\Vert u\Vert_{L^p(B_R\setminus B_{R/2})}\Vert\nabla u\Vert_{L^2({B_r})}\Vert v^*\Vert_{L^q({B_r})}+\Vert \nabla w\Vert_{L^p(B_r)}\Vert\nabla u\Vert_{L^2({B_r})}\Vert v^*\Vert_{L^q({B_r})}\right)\\
		&\leq \frac{c}{r-\rho}\Vert u\Vert_{L^p(B_R\setminus B_{R/2})}\Vert\nabla u\Vert_{L^2({B_r})}\Vert v^*\Vert_{L^q({B_r})},
	\end{aligned}
\end{equation}
with $ \frac1q+\frac1p+\frac12=1$. It follows from $ \rho u\in \dot{B}^{\frac3p-\frac32}_{\infty,\infty} $ that
\begin{equation}\label{eq2.th2.5}
\Vert v^*\Vert_{L^q({B_r})}\leq cr^{\frac3q}\left(\Vert v-v(0)\Vert_{L^{\infty}(B_r)}\right)\leq cr^{\frac3q+\frac3p-\frac12}\Vert v\Vert_{\dot{B}^{\frac3p-\frac12}_{\infty,\infty}}\leq cR\Vert \rho u\Vert_{\dot{B}^{\frac3p-\frac32}_{\infty,\infty}}.
\end{equation}
Combining \eqref{eq2.th2.3}-\eqref{eq2.th2.5} yields
\begin{equation*}
  \Phi_u(\rho) \leq \frac14\Phi_u(r)+\frac{cR^{\frac{3p-6}{p}}}{(r-\rho)^2}\Vert  u\Vert^2_{L^{p}(B_{2R}\setminus B_R)}+\frac{cR^2}{(r-\rho)^2}\Vert  u\Vert^2_{L^{p}(B_{2R}\setminus B_R)}\Vert \rho u\Vert^2_{\dot{B}^{\frac3p-\frac32}_{\infty,\infty}}.
\end{equation*}
By suitable iteration, we deduce Caccioppoli type inequality
\begin{equation*}
  \mu\int_{B_{R/2}}|\nabla u|^2\mathrm{d}x+(\mu+\lambda)\int_{B_{R/2}}|\mathrm{div}u|^2\mathrm{d}x\leq c\left(R^{\frac{p-6}{p}}+\Vert \rho u\Vert^2_{\dot{B}^{\frac3p-\frac32}_{\infty,\infty}}\right)\Vert  u\Vert^2_{L^{p}(B_{2R}\setminus B_R)}.
\end{equation*}
 Passing the limit as $ R\to\infty $, the second part of theorem is completed.

\end{proof}

\begin{proof}[{\rm \textbf{Proof of Theorem \ref{thNS3}}}]
	Let $ \phi,w $ be defined in the proof of Theorem \ref{thNS1}. Recall the function $ \phi $ has the properties:
\begin{equation*}
  	 |\nabla\phi|\leq\frac{c}{r-\rho},\quad\text{and}\quad|\nabla^2\phi|\leq\frac{c}{(r-\rho)^2}.
\end{equation*}
	 According to the interpolation and Poincar\'{e}'s inequality, one obtains
	\begin{equation}\label{eq2.th3.1}
		\begin{aligned}
			\Vert\nabla w\Vert_{L^{2q',2}(B_r)}&\leq c\Vert u\cdot\nabla\phi\Vert_{L^{2q',2}(B_r)}\leq c\Vert u\cdot\nabla\phi\Vert_{L^{2}(B_r)}^{1-\frac3{2q}}\Vert u\cdot\nabla\phi\Vert_{L^{6}(B_r)}^{\frac3{2q}}\\
			&\leq \frac{c}{(r-\rho)^{1-\frac3{2q}}}\Vert u\Vert_{L^{2}(B_r)}^{1-\frac3{2q}}\Vert\nabla (u\cdot\nabla\phi)\Vert_{L^{2}(B_r)}^{\frac3{2q}}\\
			&\leq\frac{c}{r-\rho}\Vert u\Vert_{L^{2}(B_R\setminus B_{R/2})}^{1-\frac3{2q}}\Vert\nabla u\Vert_{L^{2}(B_r)}^{\frac3{2q}}+\frac{c}{(r-\rho)^{1+\frac3{2q}}}\Vert u\Vert_{L^{2}(B_R\setminus B_{R/2})},
		\end{aligned}
	\end{equation}	
and
\begin{equation}\label{eq2.th3.2}
	\begin{aligned}
		\Vert u\Vert_{L^{2q',2}(\Omega)}&\leq c\Vert u\Vert_{L^{2}(\Omega)}^{1-\frac3{2q}}\Vert u\Vert_{L^{6}(\Omega)}^{\frac3{2q}}\\
		&\leq c\Vert u\Vert_{L^{2}(\Omega)}^{1-\frac3{2q}}\Vert\nabla u\Vert_{L^{2}(\Omega)}^{\frac3{2q}}+c\Vert u\Vert_{L^{2}(\Omega)}^{1-\frac3{2q}}\Vert u\Vert_{\dot{M}^{s,6}}^{\frac3{2q}},
	\end{aligned}
\end{equation}
with $ \Omega=B_r,B_R\text{ or }B_R\setminus B_{R/2},$ where the last inequality holds from H\"{o}lder inequality
\begin{equation*}
  \Vert (u)_{\Omega}\Vert_{L^{6}(\Omega)}\leq cR^{-\frac52}\left|\int_{B_R}u\mathrm{d}x\right|\leq CR^{\frac12-\frac3s}\Vert u\Vert_{L^s(B_R)}\leq\Vert u\Vert_{\dot{M}^{s,6}}.
\end{equation*}

Recall the energy identity \eqref{eq2.th1.3} in the proof of Theorem \ref{thNS1}
\begin{equation}\label{eq2.th3.3}
	\begin{aligned}
		\mu\int_{B_r}\phi|\nabla u|^2&\mathrm{d}x+(\mu+\lambda)\int_{B_r}\phi|\mathrm{div}u|^2\mathrm{d}x\\
		&=\frac12\int_{B_r}\rho u|u|^2 \cdot\nabla\phi\mathrm{d}x-\int_{B_r} \rho u\cdot(u\cdot\nabla w)\mathrm{d}x-\mu\int_{B_r}\nabla u:(u\otimes\nabla \phi)\mathrm{d}x\\
		&\qquad-\mu\int_{B_r}\nabla u:\nabla w\mathrm{d}x-\int_{B_r}\frac{(\mu+\lambda)}{\gamma-1}\mathrm{div}u\cdot u\cdot\nabla\phi\mathrm{d}x\\
		& = I_1+I_2+I_3+I_4+I_5.
	\end{aligned}
\end{equation}
According to H\"{o}lder's inequality of the Lorentz spaces and inequality \eqref{eq2.th3.2}, it holds that
	\begin{align*}
		|I_1|&\leq \frac{c}{r-\rho}\Vert\rho u\Vert_{L^{q,\infty}(B_r)}\Vert u\Vert^2_{L^{2q',2}(B_r\setminus B_{\rho})}\\
		&\leq\frac{cR^{\frac3q-1}}{r-\rho}\Vert \rho u\Vert_{\dot M_w^{q,3}\left( {\mathbb{R}^3} \right)}\Vert u\Vert_{L^{2}(B_R\setminus B_{R/2})}^{2-\frac3{q}}(\Vert\nabla u\Vert_{L^{2}(B_r)}^{\frac3{q}}+1).
	\end{align*}
Note that
\begin{equation}\label{eq2.th3.4}
	\Vert u\Vert_{L^{2}(B_R\setminus B_{R/2})}^{2-\frac3{q}}\leq cR^{\frac{3s-6}{2s}}\Vert u\Vert_{L^{s}(B_R\setminus B_{R/2})}^{2-\frac3{q}},
\end{equation}
which gives rise to
\begin{equation}\label{nI1}
  \begin{aligned}
	|I_1|&\leq c\Vert u\Vert_{\dot{M}^{s,6}}^{2-\frac3q}\frac{R}{r-\rho}\Vert\nabla u\Vert_{L^{2}(B_r)}^{\frac3{q}}+c\Vert u\Vert_{\dot{M}^{s,6}}^{2-\frac3q}\frac{R}{r-\rho}\\
	&\leq \frac{1}{32} \Phi_u(r)+ c\left(\frac{R}{r-\rho}\right)^{\frac{2q}{2q-3}}+c\frac{R}{r-\rho},
\end{aligned}
\end{equation}
where $ \Phi_u $ is defined in the proof of Theorem \ref{thNS1}.
Thanks to H\"{o}lder's inequality of the Lorentz spaces, it follows from \eqref{eq2.th3.1} and \eqref{eq2.th3.2} that
\begin{align*}
	|I_2|&\leq\Vert\rho u\Vert_{L^{q,\infty}(B_r)}\Vert u\Vert_{L^{2q',2}(B_r)}\Vert \nabla w\Vert_{L^{2q',2}(B_r)}\\
	&\leq \frac{cNR^{\frac3q-1}}{r-\rho}\Vert u\Vert_{L^2(B_R)}^{2-\frac3q}\Vert\nabla u\Vert_{L^{2}(B_r)}^{\frac3q}+\frac{cNR^{\frac3q-1}}{r-\rho}\Vert u\Vert_{L^2(B_R)}^{2-\frac3q}\Vert\nabla u\Vert_{L^{2}(B_r)}^{\frac3{2q}}\\
	&\qquad+\frac{cNR^{\frac3q-1}}{(r-\rho)^{1+\frac3{2q}}}\Vert u\Vert_{L^2(B_R)}^{2-\frac3{2q}}\Vert\nabla u\Vert_{L^{2}(B_r)}^{\frac3{2q}}+\frac{cNR^{\frac3q-1}}{(r-\rho)^{1+\frac3{2q}}}\Vert u\Vert_{L^2(B_R)}^{2-\frac3{2q}}\\
	&=I_2^1+I_2^2+I_2^3+I_2^4.
\end{align*}
 After a direct calculation, one obtained
\begin{align*}
	I_2^i&\leq\frac{1}{32}\Phi_u(r)+c\left(\frac{R}{r-\rho}\right)^{\kappa_i},\quad\text{for}\quad i=1,2,3,4,
\end{align*}
with $\kappa_1=\frac{2q}{2q-3},\, \kappa_2=\frac{4q}{4q-3},\, \kappa_3=\frac{4q+6}{4q-3}$\, and\, $\kappa_4=\frac{2q+3}{2q}$.
Summmarying above inequalities yields
\begin{align}\label{nI2}
	|I_2|\leq\frac18\Phi_u(r)+c\left(\left(\frac{R}{r-\rho}\right)^{\kappa_1}+\left(\frac{R}{r-\rho}\right)^{\kappa_2}+\left(\frac{R}{r-\rho}\right)^{\kappa_3}
+\left(\frac{R}{r-\rho}\right)^{\kappa_4}\right).
\end{align}
Thanks to the property of $w$, it holds that
\begin{equation}\label{nI3}
  |I_3|\leq \frac{cR^{\frac{3s-6}{2s}}}{r-\rho}\Vert u\Vert_{L^{s}(B_R)}\Vert\nabla u\Vert_{L^{2}(B_r)}\leq\frac{1}{32}\Phi_u(r)+c\left(\frac{R}{r-\rho}\right)^2.
\end{equation}
and
\begin{equation}\label{nI4}
  \begin{aligned}
	|I_4|+|I_5|&\leq cR^{\frac{3s-6}{2s}}\Vert\nabla w\Vert_{L^{s}(B_r)}\Vert\nabla u\Vert_{L^{2}(B_r)}+\frac{CR^{\frac{3s-6}{2s}}}{r-\rho}\Vert u\Vert_{L^{s}(B_R)}\Vert\mathrm{div} u\Vert_{L^{2}(B_r)}\\
	&\leq\frac{cR^{\frac{3s-6}{2s}}}{r-\rho}\Vert u\Vert_{L^{s}(B_R)}\Vert\nabla u\Vert_{L^{2}(B_r)}\frac{CR^{\frac{3s-6}{2s}}}{r-\rho}\Vert u\Vert_{L^{s}(B_R)}\Vert\mathrm{div} u\Vert_{L^{2}(B_r)}\\
&\leq\frac{1}{32}\Phi_u(r)+c\left(\frac{R}{r-\rho}\right)^2.
\end{aligned}
\end{equation}
Combining \eqref{eq2.th3.3} and \eqref{nI1}-\eqref{nI4}, one obtains
\begin{equation*}
  \Phi_u(\rho)\leq\frac14\Phi_u(r)+cf\left(\frac{R}{r-\rho}\right),\quad \frac{R}2\leq\rho<r\leq R,
\end{equation*}
with $c=c(s,q,\Vert u\Vert_{\dot{M}^{s,6}}, \Vert\rho u\Vert_{\dot M_w^{p,3}})$\, and \, $f(z)=z^{\kappa_1}+z^{\kappa_2}+z^{\kappa_3}+z^{\kappa_4}+z^2+z$. Thus, the iteration lemma implies Caccioppoli type inequality
\begin{equation*}
  \int_{B_{R/2}}(\mu|\nabla u|^2+(\mu+\lambda)|\mathrm{div}u|^2)\mathrm{d}x\leq c(s,q,\Vert\rho u\Vert_{\dot M_w^{p,3}},\Vert u\Vert_{\dot{M}^{s,6}}),
\end{equation*}
which implies $u\in \dot{H}^1(\mathbb{R}^3)$.

Now the Sobolev embedding shows that there is a constant vector $C$ such that $u-C \in L^6(\mathbb{R}^3)$. We claim that $C=0.$ In fact,
if $C\ne 0$, it is obvious that
\begin{equation*}
\begin{aligned}
 {R^{3\left( {\frac{1}{6} - \frac{1}{s}} \right)}}{\left\| u \right\|_{{L^s}({B_R}(x))}} &\ge {R^{3\left( {\frac{1}{6} - \frac{1}{s}} \right)}}{\left\| C \right\|_{{L^s}({B_R}(x))}} - {R^{3\left( {\frac{1}{6} - \frac{1}{s}} \right)}}{\left\| {u - C} \right\|_{{L^s}({B_R}(x))}}\\
 &\ge c{R^2} - {\left\| {u - C} \right\|_{{L^6}(\mathbb{R}^3)}}\to\infty,\quad\text{as}\quad R\to\infty.
 \end{aligned}
\end{equation*}
It is a contradiction since $u\in \dot{M}^{s,6}(\mathbb{R}^3)$.

Taking $r=R$ and $\rho =\frac{R}{2}$ in \eqref{eq2.th3.3} yields
 \begin{equation}\label{eq2.th3.5}
	\begin{aligned}
		\mu\int_{B_R}\phi|\nabla u|^2&\mathrm{d}x+(\mu+\lambda)\int_{B_R}\phi|\mathrm{div}u|^2\mathrm{d}x
		= I_1(R)+I_2(R)+I_3(R)+I_4(R)+I_5(R).
	\end{aligned}
\end{equation}
Due to H\"{o}lder's inequality, it holds that
	\begin{equation}\label{eq2.th3.6}
		\begin{aligned}
		|I_1(R)|&\leq cR^{-1}\Vert\rho u\Vert_{L^{q,\infty}(B_R)}\Vert u\Vert^2_{L^{2q',2}(B_R\setminus B_{R/2})}\\
		&\leq {cR^{\frac3q-2}}\Vert u\Vert_{L^{2}(B_R\setminus B_{R/2})}^{2-\frac3{q}}(\Vert\nabla u\Vert_{L^{2}(B_R\setminus B_{R/2})}^{\frac3{q}}+1) \\
       &\leq c\Vert u\Vert_{L^{6}(B_R\setminus B_{R/2})}^{2-\frac3{q}}(\Vert\nabla u\Vert_{L^{2}(B_R\setminus B_{R/2})}^{\frac3{q}}+1).
	\end{aligned}
	\end{equation}
Using the same argument yields
\begin{equation}\label{eq2.th3.7}
\begin{aligned}
  |I_2(R)&|\leq c\Big(\Vert u\Vert_{L^6(B_R\setminus B_{R/2})}^{1-\frac3{2q}}\Vert u\Vert_{L^6(B_R)}^{1-\frac3{2q}}\Vert\nabla u\Vert_{L^{2}(B_R)}^{\frac3q}+\Vert u\Vert_{L^6(B_R\setminus B_{R/2})}\Vert u\Vert_{L^6(B_R)}^{1-\frac3{2q}}\\
 	&\qquad+\Vert u\Vert_{L^6(B_R\setminus B_{R/2})}^{1-\frac3{2q}}\Vert u\Vert_{L^6(B_R)}^{1-\frac3{2q}}\Vert\nabla u\Vert_{L^{2}(B_R)}^{\frac3{2q}}\\
    &\quad\quad+\Vert u\Vert_{L^6(B_R\setminus B_{R/2})}\Vert u\Vert_{L^6(B_R)}^{1-\frac3{2q}}\Vert\nabla  u\Vert_{L^{2}(B_R)}^{\frac3{2q}}\Big),
 \end{aligned}
\end{equation}
and
 \begin{align}\label{eq2.th3.8}
 	|I_3(R)|+|I_5(R)|\leq c\Vert u\Vert_{L^6(B_R\setminus B_{R/2})}\Vert\nabla u\Vert_{L^2(B_R)}.
 \end{align}
It follows from \eqref{eq2.th1.1} that
\begin{equation}\label{eq2.th3.9}
  |I_4(R)|\leq cR\Vert \nabla w\Vert_{L^6(B_R)}\Vert\nabla u\Vert_{L^2(B_R)}\leq c\Vert u\Vert_{L^6(B_R\setminus B_{R/2})}\Vert\nabla u\Vert_{L^2(B_R)}.
\end{equation}
Applying inequalities \eqref{eq2.th3.6}-\eqref{eq2.th3.9} gives $\mathop {\lim }_{R \to  + \infty } {I_k}\left( R \right) = 0$ for $k=1,2,3,4,5$, which implies
\begin{equation*}
\begin{aligned}
  \mu &\int_{\mathbb{R}^3} |\nabla u{|^2}dx + (\mu  + \lambda )\int_{\mathbb{R}^3}  |divu{|^2}dx\\
  &\quad\quad=\mathop {\lim }\limits_{R \to \infty } (\mu\int_{B_R}\phi|\nabla u|^2\mathrm{d}x+(\mu+\lambda)\int_{B_R}\phi|\mathrm{div}u|^2\mathrm{d}x)=0.
   \end{aligned}
\end{equation*}
Since $u \in L^6({\mathbb{R}^3})$, the proof of Theorem \ref{thNS3} is completed.
\end{proof}


\section*{Declaration of competing interest}
We thank professor Hailiang Li for many helpful discussions and continual interest in this work. Li's research is partially supported by the National Natural Science Foundation of China grants (No. 11931010, No. 12226326), and the key research project of Academy for Multidisciplinary Studies, Capital Normal University. Quansen Jiu was partially supported by National Natural Science Foundation
of China under grants (No. 11931010, No. 12061003).

\end{document}